# A strong uniform convergence rate of a kernel conditional quantile estimator under random left-truncation and dependent data

Elias Ould-Saïd[1], Djabrane Yahia[2], Abdelhakim Necir[2]

[1] *Univ. Lille Nord de France, L.M.P.A. J. Liouville (Univ. du Littoral Côte d'Opale), BP 699, 62228 Calais, France.*
*e-mail:* `ouldsaid@lmpa.univ-littoral.fr`*(corresponding author)*

[2] *Lab. of Applied Mathematics, Univ. Mohamed Khider, PB 145, 07000 Biskra, Algérie.*
(`yahia_dj` or `necirabdelhakim`)`@yahoo.fr`

**Abstract:** In this paper we study some asymptotic properties of the kernel conditional quantile estimator with randomly left-truncated data which exhibit some kind of dependence. We extend the result obtained by Lemdani, Ould-Saïd and Poulin [16] in the iid case. The uniform strong convergence rate of the estimator under strong mixing hypothesis is obtained.

**Keywords and phrases:** Kernel estimator · Quantile function· Rate of convergence · Strong mixing · Strong uniform consistency· Truncated data.

**Mathematics Subject classification (2000)** · 62G05 · 62G20.

## Contents

## 1. Introduction

Let $\mathcal{Y}$ and $\mathcal{T}$ be two real random variables (rv) with unknown cumulative distribution functions (df) $F$ and $G$ respectively, both assumed to be contin-uous. Let $\mathcal{X}$ be a real-valued random covariable with df $V$ and continuous density $v$. Under random left-truncation (RLT), the rv of interest $\mathcal{Y}$ is interfered by the truncation rv $\mathcal{T}$, in such a way that $\mathcal{Y}$ and $\mathcal{T}$ are observed only if $\mathcal{Y} \geq \mathcal{T}$. Such data occur in astronomy and economics (see Woodroofe [31], Feigelson and Babu [7], Wang *et al.* [30], Tsai *et al.* [29]) and also in epidemiology and biometry (see, e.g., He and Yang [12]).

If there were no truncation, we could think of the observations as $(\mathcal{X}_j, \mathcal{Y}_j, \mathcal{T}_j)$; $1 \leq j \leq N$, where the sample size $N$ is deterministic, but unknown. Under RLT, however, some of these vectors would be missing and for notational convenience, we shall denote $(X_i, Y_i, T_i)$; $1 \leq i \leq n$, $(n \leq N)$ the observed subsequence subject to $Y_i \geq T_i$ from the $N-$sample. As a consequence of truncation, the size of actually observed sample, $n$, is a binomial rv with parameters $N$ and $\mu := \mathbb{P}(\mathcal{Y} \geq \mathcal{T}) > 0$. By the strong law of large numbers we have, as $N \to \infty$

$$\mu_n := \frac{n}{N} \to \mu, \quad \mathbb{P}-a.s. \tag{1}$$

Now we consider the joint df $F(.,.)$ of the random vector $(\mathcal{X}, \mathcal{Y})$ related to the $N-$sample and suppose it is of class $\mathcal{C}^1$. The conditional df of $\mathcal{Y}$ given $\mathcal{X} = x =: (x_1, ..., x_d)^t$, that is $F(y|x) = \mathbb{E}\left[\mathbf{1}_{\{\mathcal{Y} \leq y\}} | \mathcal{X} = x\right]$ which may be rewritten into

$$F(.|x) = \frac{F_1(x,.)}{v(x)} \tag{2}$$

where $F_1(x,.)$ is the first derivative of $F(x,\cdot)$ with respect to $x$. For all fixed $p \in (0, 1)$, the $p^{th}$ conditional quantile of $F$ given $\mathcal{X} = x$ is defined by

$$q_p(x) := \inf\{y \in \mathbb{R} : F(y|x) \geq p\}.$$

It is well known that the quantile function can give a good description of the data (see, Chaudhuri *et al.* [5]), such as robustness to heavy-tailed error distributions and outliers, especially the conditional median function $q_{1/2}(x)$ for asymmetric distribution, which can provide a useful alternative to the ordinary regression based on the mean. The nonparametric estimation of conditional quantile has first been considered in the case of complete data (no truncation). Roussas [24] showed the convergence and asymptotic normality of kernel estimates of conditional quantile under Markov assumptions. For independent and identically distributed (iid) rv's, Stone [27] proved the weak consistency of kernel estimates. The uniform consistency was studied by Schlee [26] and Gannoun [9]. The asymptotic normality has been established by Samanta [25]. Mehra *et al.* [20] proposed and discussed certain smooth variants (based both on single as well as double kernel weights) of the standard conditional quantile estimator, proved the asymptotic normality and found an almost sure (a.s.) convergence rate, whereas Xiang [32] gave the asymptotic normality and a law of the iterated logarithm for a new kernel estimator. In the dependent case, the convergence of nonparametric estimation of quantile was proved by Gannoun [10] and Boente and Fraiman [1].

In the RLT model, Gürler *et al.* [11] gave a Bahadur-type representation for the quantile function and asymptotic normality. Its extension to time series analysis was obtained by Lemdani *et al.* [15].

The aim of this paper is to establish a strong uniform convergence rate for the kernel conditional quantile estimator with randomly left-truncated data under $\alpha-$mixing conditions whose definition is given below. Hence, we extend the obtained result by Lemdani *et al.* [16] in the iid case.

First, let $\mathcal{F}_i^k(Z)$ denotes the $\sigma$-field of events generated by $\{Z_j,\ i \leq j \leq k\}$. For easy reference, let us recall the following definition.

**Definition 1.1** *Let $\{Z_i,\ i \geq 1\}$ denotes a sequence of random variables. Given a positive integer $n$, set:*

$$\alpha(n) = \sup\left\{|\mathbb{P}(A \cap B) - \mathbb{P}(A)\mathbb{P}(B)| : A \in \mathcal{F}_1^k(Z),\ B \in \mathcal{F}_{k+n}^\infty(Z),\ k \in \mathbb{N}\right\}.$$

*The sequence is said to be $\alpha-$mixing (strongly mixing) if the mixing coefficient $\alpha(n) \to 0$.*

Among various mixing conditions used in the literature, $\alpha-$mixing is reasonably weak and has many practical applications (see, e.g. Doukhan [6] or Cai ([3, 4] for more details). In particular, Masry and Tj$\phi$stheim [18] proved that, both ARCH processes and nonlinear additive AR models with exogenous variables, which are particularly popular in finance and econometrics, are stationary and $\alpha-$mixing.

The rest of the paper is organized as follows. In Section 2, we recall a definition of the kernel conditional quantile estimator with randomly left-truncated data. Assumptions and main results are given in Section 3. Section 4 is devoted to application to prediction. Finally, the proofs of the main results are postponed to Section 5 with some auxiliary results and their proofs.

## 2. Definition of the estimator

In the sequel, the letters $C$ and $C'$ are used indiscriminately as generic constants. Note also that, $N$ is unknown and $n$ is known (although random), our results will not be stated with respect to the probability measure $\mathbb{P}$ (related to the $N-$sample) but will involve the conditional probability $\mathbf{P}$ (related to the $n-$sample). Also $\mathbb{E}$ and $\mathbf{E}$ will denote the expectation operators related to $\mathbb{P}$ and $\mathbf{P}$, respectively. Finally, we denote by a superscript $(^*)$ any df that is associated to the observed sample.

The estimation of conditional df is based on the choice of weights. For the complete data, the well-known Nadaraya-Watson weights are given by

$$W_{i,N}(x) = \frac{K\left\{(x-X_i)/h_N\right\}}{\sum_{i=1}^N K\left\{(x-X_i)/h_N\right\}} = \frac{(Nh_N)^{-1} K\left\{(x-x_i)/h_N\right\}}{v_N(x)} \quad (3)$$

that are measurable functions of $x$ depending on $X_1, ..., X_N$, with the convention $0/0 = 0$. The kernel $K$ is a measurable function on $\mathbb{R}^d$ and $(h_N)$ a nonnegative sequence which tends to zero as $N$ tends to infinity. The regression estimator based on the $N$-sample is then given by

$$r_N(x) = \frac{(Nh_N)^{-1} \sum_{i=1}^n Y_i K\left\{(x-X_i)/h_N\right\}}{v_N(x)} \quad (4)$$

where $v_N$ is a well knwon kernel estimator of $v$ based on the $N-$sample. As $N$ is unknown, then $v_N(\cdot)$ cannot be calculated and therefore $r_N(\cdot)$. On the other hand, based on the $n$-sample, the kernel estimator

$$v_n^*(x) = \frac{1}{nh_n} \sum_{i=1}^n K\left(\frac{x-X_i}{h_n}\right) \quad (5)$$

is an estimator of the conditional density $v^*(x)$ (given $\mathcal{Y} \geq \mathcal{T}$), see Ould-Saïd and Lemdani [21].

Under RLT sampling scheme, the conditional joint distribution (Stute, [28]) of $(Y, T)$ becomes

$$J^*(y,t) = \mathbf{P}(Y \leq y, T \leq t) = \mathbb{P}(Y \leq y, T \leq t | Y \geq T)$$
$$= \mu^{-1} \int_{-\infty}^{y} G(t \wedge u) dF(u)$$

where $t \wedge u := \min(t, u)$. The marginal distribution and their empirical versions are defined by

$$F^*(y) = \mu^{-1} \int_{-\infty}^{y} G(u) dF(u), \qquad F_n^*(y) = n^{-1} \sum_{i=1}^{n} \mathbf{1}_{\{Y_i \leq y\}},$$

$$G^*(t) = \mu^{-1} \int_{-\infty}^{\infty} G(t \wedge u) dF(u) \quad \text{and} \quad G_n^*(t) = n^{-1} \sum_{i=1}^{n} \mathbf{1}_{\{T_i \leq t\}},$$

where $\mathbf{1}_A$ denote the indicator function of the set $A$.
In the sequel we use the following consistent estimator

$$\mu_n = \frac{G_n(y) \left[ (1 - F_n(y-)) \right]}{C_n(y)}, \tag{6}$$

for any $y$ such that $C_n(y) \neq 0$, where $F_n(y-)$ denotes the left-limite of $F_n$ at $y$. Here $F_n$ and $G_n$ are the product-limit estimators (Lynden-Bell [17]) for $F$ and $G$, respectively i.e.,

$$F_n(y) = 1 - \prod_{i/Y_i \leq y} \left[ \frac{nC_n(Y_i) - 1}{nC_n(Y_i)} \right], \quad G_n(y) = \prod_{i/T_i > y} \left[ \frac{nC_n(T_i) - 1}{nC_n(T_i)} \right],$$

where $C_n(y) = n^{-1} \sum_{i=1}^{n} \mathbf{1}_{\{T_i \leq y \leq Y_i\}}$ is the empirical estimator of

$$C(y) = \mathbb{P}(T \leq y \leq Y | Y \geq T).$$

He and Yang [13] proved that $\mu_n$ does not depend on $y$ and its value can then be obtained for any $y$ such that $C_n(y) \neq 0$. Furthermore, they showed (see their Corollary 2.5) its $\mathbb{P}-a.s.$ consistency.

Suppose now that one observes the $n$ triplets $(X_i, Y_i, T_i)$ among the $N$ ones and for any df $L$, denotes the left and right endpoint of its support by $a_L := \inf\{x : L(x) > 0\}$ and $b_L := \sup\{x : L(x) < 1\}$, respectively. Then under the current model, as discussed by Woodroofe [31], $F$ and $G$ can be estimated completely only if

$$a_G \leq a_F, \quad b_G \leq b_F \quad \text{and} \quad \int_{a_F}^{\infty} \frac{dF}{G} < \infty.$$

In order to estimate the marginal density $v$ we have to take into account the truncation and the estimator

$$v_n(x) = \frac{\mu_n}{n h_n} \sum_{i=1}^{n} \frac{1}{G_n(Y_i)} K\left(\frac{x - x_i}{h_n}\right) \tag{7}$$

is considered in Ould-Saïd and Lemdani [15]. Note that in this formula and the forthcoming, the sum is taken only for $i$ such that $G_n(Y_i) \neq 0$.

Then, adapting Ould-Saïd-Lemdani's weights, we get the following estimator of the conditional df of $\mathcal{Y}$ given $\boldsymbol{\mathcal{X}} = x$

$$F_n(y|x) = \mu_n \sum_{i=1}^n \widetilde{W}_{i,n}(x) G_n^{-1}(Y_i) H\left(\frac{y-Y_i}{h_n}\right)$$

$$= \frac{\sum_{i=1}^n \frac{1}{G_n(Y_i)} K\left(\frac{x-X_i}{h_n}\right) H\left(\frac{y-Y_i}{h_n}\right)}{\sum_{i=1}^n \frac{1}{G_n(Y_i)} K\left(\frac{x-X_i}{h_n}\right)}$$

$$=: \frac{F_{1,n}(x,y)}{v_n(x)} \qquad (8)$$

where $H$ is a distribution function defined on $\mathbb{R}$, and

$$F_{1,n}(x,y) = \frac{\mu_n}{nh_n^d} \sum_{i=1}^n \frac{1}{G_n(Y_i)} K\left(\frac{x-x_i}{h_n}\right) H\left(\frac{y-Y_i}{h_n}\right) \qquad (9)$$

is an estimator of $F_1(x,y)$. As the latter is continuous, it is clear that it is better to define a smooth estimator by using a continuous function $H(\cdot)$ instead of a step function $I_{\{\cdot\}}$. We point out here that the estimators (8) and (9) have been already defined in Lemdani *et al.* [16].

Then a natural estimator of the $p^{th}$ conditional quantile $q_p(x)$ is given by

$$q_{p,n}(x) := \inf\{y \in \mathbb{R} : F_n(y|x) \geq p\}. \qquad (10)$$

which satisfies $F_n(q_{p,n}(x)|x) = p$.

## 3. Assumptions and main results

In what follows, we focus our attention on the case of a univariate covariable ($d = 1$) and denote $X$ for $x$ and $K$ for $K_1$. Assume that $0 = a_G < a_F$ and $b_G \leq b_F$. We consider two real numbers $a$ and $b$ such that $a_F < a < b < b_F$. Let $\Omega$ be a compact subset of $\Omega_0 = \{x \in \mathbb{R} | v(x) > 0\}$ and $\gamma := \inf_{x \in \Omega} v(x) > 0$.

We introduce some assumptions, gathered below for easy reference needed to state our results.

(K1) $K$ is a positive-valued, bounded probability density, Hölder continuous with exponent $\beta > 0$ and satisfying
$$|u| K(u) \to 0 \quad as \quad |u| \to +\infty.$$

(K2) $H$ is a df with $\mathcal{C}^1$-probability density $H^{(1)}$ which is poisitive, bounded and has compact support. It is also Hölderian with exponent $\beta$.

(K3) i) $H^{(1)}$ and $K$ are second-order kernels,
ii) $\int K^2(r)dr < \infty$.

(M1) $\{(X_i, Y_i); i \geq 1\}$ is a sequence of stationary $\alpha$-mixing random variables with coefficient $\alpha(n)$.

(M2) $\{T_i; i \geq 1\}$ is a sequence of iid truncating variables independent of $\{(X_i, Y_i), i \geq 1\}$ with common continuous df $G$.

(M3) There exists $\nu > 5 + 1/\beta$ for some $\beta > 1/7$ such that
$\forall n, \alpha(n) = O(n^{-\nu})$.

(D1) The conditional density $v^*(.)$ is twice continuously differentiable.

(D2) The joint conditional density $v^*(.,.)$ of $(X_i, X_j)$ exists and satisfies
$$\sup_{r,s} |v^*(r,s) - v^*(r)v^*(s)| \leq C < \infty,$$
for some constant $C$ not depending on $(i, j)$.

(D3) The joint conditional density of $(X_i, Y_i, X_j, Y_j)$ denoted by $f^*(.,.,.,.)$, exists and satisfies for any constant $C$,
$$\sup_{r,s,t,u} |f^*(r,s,t,u) - f^*(r,s)f^*(t,u)| \leq C < \infty.$$

(D4) The joint density $f(.,.)$ is bounded and twice continuously differentiable.

(D5) The marginal density $v(.)$ is locally Lipschitz continuous over $\Omega_0$.

The bandwidth $h_n =: h$ satisfies:

(H1)
$$h \downarrow 0, \quad \frac{\log n}{nh} \to 0 \quad \text{and} \quad h = o(1/\log n), \qquad \text{as } n \to \infty,$$

(H2)
$$Cn^{\frac{(3-\nu)\beta}{\beta(\nu+1)+4\beta+1}+\eta} < h < C'n^{\frac{1}{1-\nu}},$$
where $\eta$ satisfies
$$\frac{2}{\beta(\nu+1)+4\beta+1} < \eta < \frac{(\nu-3)\beta}{\beta(\nu+1)+4\beta+1} + \frac{1}{1-\nu},$$

$\nu$ and $\beta$ are the same as in $(M3)$.

**Remark 3.1** *Assumptions $(K)$ are quite usual in kernel estimation. Conditions $(D1)$, $(D4)$ and $(D5)$ are needed in the study of the bias term. $(D2)$ and $(D3)$ are needed for covariance calculus and take similar forms to those used under mixing. Hypothesis $(H2)$ is used in Ould-Saïd and Tatachak [22] and is needed to establish Lemma 5.1 and Lemma 5.4. Assumptions $(M)$ concern the mixing processes structure which are standard in such situation. The choice of $\beta$ seems to be surprising, but it is only technical choice which permit us to make one of the variance term to be negligible.*

**Remark 3.2** *Here we point out that we can not suppose that the original data (that is the $N$-sample) satisfies some kind of dependency. Indeed, we do not know if the observed data are $\alpha$-mixing or are not. And if they are, we do not know the coefficient. Therefore, we suppose that the observed data satisfy some kind of mixing condition.*

**Remark 3.3** *As we are interested in the number $n$ of observations ($N$ is unknown), we give asymptotics as $n \to \infty$ unless otherwise specified. Since $n \leq N$, this implies $N \to \infty$ and these results also hold under $\mathbb{P} - a.s.$ as $N \to \infty$.*

Our first result, stated in Proposition 3.1, is the uniform almost sure convergence with rate of the conditional df estimator defined in (8).

**Proposition 3.1** *Under assumptions $(K)$, $(M)$, $(D)$ and $(H)$, we have*

$$\sup_{x \in \Omega} \sup_{a \leq y \leq b} |F_n(y|x) - F(y|x)| = O\left(\max\left\{\sqrt{\frac{\log n}{nh}}, h^2\right\}\right), \quad \mathbf{P}-a.s. \quad as \quad n \to \infty$$

The second result deals with the strong uniform convergence with rate of the kernel conditional quantile estimator $q_{p,n}(.)$ which is given in the following theorem.

**Theorem 3.1** *Under the assumptions of Proposition 3.1 and for each fixed $p \in (0,1)$ if the function $q_p$ satisfies for given $\varepsilon > 0$ there exists $\beta > 0$ such that*

$$\forall \eta_p : \Omega \to \mathbb{R}, \ \sup_{x \in \Omega} |q_p(x) - \eta_p(x)| \geq \varepsilon \Rightarrow \sup_{x \in \Omega} |F(q_p(x)) - F(\eta_p(x))| \geq \beta, \tag{11}$$

we have
$$\lim_{n \to \infty} \sup_{x \in \Omega} |q_{p,n}(x) - q_p(x)| = 0, \quad \mathbf{P} - a.s.$$

*Furthermore, we have*

$$\sup_{x \in \Omega} |q_{p,n}(x) - q_p(x)| = O\left(\max\left\{\sqrt{\frac{\log n}{nh}}, h^2\right\}\right), \quad \mathbf{P} - a.s. \quad as \quad n \to \infty$$

## 4. Applications to prediction

It is well known, from the robustness theory that the median is more robust than the mean, therefore the conditional median, $\mu(x) = q_{1/2}(x)$, is a good alternative to the conditional mean as a predictor for a variable $Y$ given $X = x$. Note that the estimation of $\mu(x)$ is given by $\mu_n(x) = q_{\frac{1}{2},n}(x)$. Using this considerations and section 2, we want to predict the non observed r.v. $Y_{n+1}$ (which corresponds to some modality of our problem), from available data $X_1, \ldots, X_n$. Given a new value $X_{n+1}$, we can predict the corresponding response $Y_{n+1}$ by
$$\widehat{Y}_{n+1} = \mu_n(X_{n+1}) = q_{1/2,n}(X_{n+1}).$$

Nevertheless to say, that the theoretical predictor is given by $\mu(X_{n+1}) = q_{1/2}(X_{n+1})$.

Applying the above Theorem, we have the following corollary:

**Corollary 4.1** *Under the assumptions of Theorem 3.1, we have*

$$\left|q_{1/2,n}(X_{n+1}) - q_{1/2}(X_{n+1})\right| \longrightarrow 0, \quad \mathbf{P} - a.s. \quad as \quad n \to \infty.$$

## 5. Proofs

We need some auxiliary results and notations to prove our results. The first lemma gives the uniform convergence with rate of the estimator $v_n^*(x)$ defined in (5).

**Lemma 5.1** *Under $(K1)$, $(K3)$, $(M)$, $(D1)$, $(D2)$ and $(H)$ we have*

$$\sup_{x \in \Omega} |v_n^*(x) - v^*(x)| = O\left(\max\left\{\sqrt{\frac{\log n}{nh}}, h^2\right\}\right), \quad \mathbf{P} - a.s. \quad as \quad n \to \infty$$

**Proof.** We have

$$\sup_{x \in \Omega} |v_n^*(x) - v^*(x)| \leq \sup_{x \in \Omega} |v_n^*(x) - \mathbf{E}[v_n^*(x)]| + \sup_{x \in \Omega} |\mathbf{E}[v_n^*(x)] - v^*(x)|$$
$$=: \mathcal{I}_{1n} + \mathcal{I}_{2n}. \tag{12}$$

We begin by study the variability term $\mathcal{I}_{1n}$. The idea consists in using an exponential inequality taking into account the $\alpha$-mixing structure. The compact set $\Omega$ can be covered by a finite number $l_n$ of intervals of length $\omega_n = (n^{-1}h^{1+2\beta})^{\frac{1}{2\beta}}$, where $\beta$ is the Hölder exponent. Let $I_k := I(x_k, \omega_n)$; $k = 1, ..., l_n$, denote each interval centered at some points $x_k$. Since $\Omega$ is bounded, there exists a constant $C$ such that $\omega_n l_n \leq C$. For any $x$ in $\Omega$, there exists $I_k$ which contains $x$ such that $|x - x_k| \leq \omega_n$. We start by writing

$$\triangle_i(x) := \frac{1}{nh} \left\{ K\left(\frac{x - X_i}{h}\right) - \mathbf{E}\left[K\left(\frac{x - X_1}{h}\right)\right] \right\}.$$

Clearly, we have

$$\sum_{i=1}^{n} \triangle_i(x) = \{(v_n^*(x) - v_n^*(x_k)) - (\mathbf{E}[v_n^*(x)] - \mathbf{E}[v_n^*(x_k)])\} + (v_n^*(x_k) - \mathbf{E}[v_n^*(x_k)])$$
$$=: \sum_{i=1}^{n} \widetilde{\triangle}_i(x) + \sum_{i=1}^{n} \triangle_i(x_k).$$

Hence

$$\sup_{x \in \Omega} \left|\sum_{i=1}^{n} \triangle_i(x)\right| \leq \max_{1 \leq k \leq l_n} \sup_{x \in I_k} \left|\sum_{i=1}^{n} \widetilde{\triangle}_i(x)\right| + \max_{1 \leq k \leq l_n} \left|\sum_{i=1}^{n} \triangle_i(x_k)\right|$$
$$=: S_{1n} + S_{2n}. \tag{13}$$

Firstly, we have under assumption $(K1)$,

$$\sup_{x \in I_k} \left|\sum_{i=1}^{n} \widetilde{\triangle}_i(x)\right| \leq \frac{1}{nh} \sum_{i=1}^{n} \left|K\left(\frac{x - X_i}{h}\right) - K\left(\frac{x_k - X_i}{h}\right)\right|$$
$$+ \frac{1}{h} \mathbf{E}\left[\left|K\left(\frac{x - X_1}{h}\right) - K\left(\frac{x_k - X_1}{h}\right)\right|\right]$$
$$\leq \frac{2 \sup_{x \in I_k} |x - x_k|^\beta}{h^{1+\beta}}$$
$$\leq C\omega_n^\beta h^{-1-\beta} = O\left((nh)^{-1/2}\right). \tag{14}$$

Hence, by $(H1)$ and for $n$ large enough, we get $S_{1n} = o_{\mathbf{P}}(1)$.

We now turn to the term $S_{2n}$ in (13). Under $(K1)$, the rv's $U_i = nh\triangle_i(x_k)$ are centered and bounded. The use of the well known Fuk-Nagaev's inequality (see Rio [23, formula 6.19b, page 87]) slightly modified in Ferraty and Vieu [8, see proposition A.11-ii), page 237], allows one to get, for all $\varepsilon > 0$ and $r > 1$

$$\mathbf{P}\left\{\max_{1\leq k\leq l_n}\left|\sum_{i=1}^{n}\triangle_i(x_k)\right| > \varepsilon\right\} \leq \sum_{k=1}^{l_n}\mathbf{P}\left\{\left|\sum_{i=1}^{n}U_i(x_k)\right| > nh\varepsilon\right\}$$

$$\leq C\omega_n^{-1}\left\{\frac{n}{r}\left(\frac{r}{\varepsilon nh}\right)^{\nu+1} + \left(1+\frac{\varepsilon^2 n^2 h^2}{rs_n^2}\right)^{-\frac{r}{2}}\right\}$$

$$=: Q_{1n} + Q_{2n} \qquad (15)$$

where

$$s_n^2 = \sum_{1\leq i\leq n}\sum_{1\leq j\leq n}|Cov(U_i, U_j)|.$$

Putting

$$r = (\log n)^{1+\delta}, \text{ where } \delta > 0, \text{ and } \varepsilon = \varepsilon_0\sqrt{\frac{\log n}{nh}}, \text{ for some } \varepsilon_0 > 0. \quad (16)$$

We have

$$Q_{1n} = C(n^{-1}h^{1+2\beta})^{\frac{-1}{2\beta}}\frac{n}{(\log n)^{1+\delta}}\left(\frac{(\log n)^{1+\delta}}{\varepsilon_0\sqrt{nh\log n}}\right)^{\nu+1}$$

$$= Cn^{1-\frac{\nu+1}{2}+\frac{1}{2\beta}}h^{-\left(\frac{1}{2\beta}+1+\frac{\nu+1}{2}\right)}(\log n)^{\nu(1+\delta)-\frac{\nu+1}{2}}\varepsilon_0^{-(\nu+1)}.$$

Note that under $(M3)$, it is easy to see that the following modified assumption $(H'2)$ of $(H2)$ hold,

$$Cn^{\frac{(3-\nu)\beta}{\beta(\nu+1)+2\beta+1}+\eta} < h < C'n^{\frac{1}{1-\nu}}, \qquad (H'2)$$

where $\eta$ satisfies

$$\frac{1}{\beta(\nu+1)+2\beta+1} < \eta < \frac{(\nu-3)\beta}{\beta(\nu+1)+2\beta+1} + \frac{1}{1-\nu}, \qquad (17)$$

$\nu$ and $\beta$ are the same as in $(M3)$.

Then, from the left-hand side of $(H'2)$

$$Q_{1n} \le C'(\log n)^{\nu(1+\delta)-\frac{\nu+1}{2}} n^{-1-\frac{\eta}{2\beta}(\beta(\nu+1)+2\beta+1-\frac{1}{\eta})}.$$

Hence, for any $\eta$ as in (17), $Q_{1n}$ is bounded by the term of a finite-sum series. Before we focus on $Q_{2n}$, we have to study the asymptotic behavior of

$$s_n^2 = \sum_{i=1}^n Var(U_i) + \sum_{i \ne j} |Cov(U_i, U_j)|$$
$$=: s_n^{var} + s_n^{cov}.$$

First, by $(K3:ii)$, $(D1)$ and a change of variable, we obtain

$$s_n^{var} = nVar(U_1)$$
$$= n \left\{ \mathbf{E}\left[K^2\left(\frac{x_k - X_1}{h}\right)\right] - \mathbf{E}^2\left[K\left(\frac{x_k - X_1}{h}\right)\right] \right\}$$
$$= O(nh). \tag{18}$$

On the other hand, a change of variable, $(K1)$, $(M1)$ and $(D2)$ lead to

$$|Cov(U_i, U_j)| = |\mathbf{E}[U_i U_j]|$$
$$\le \iint K\left(\frac{x_k - r}{h}\right) K\left(\frac{x_k - s}{h}\right) |v^*(r,s) - v^*(r)v^*(s)|\, drds$$
$$= O(h^2). \tag{19}$$

Note also that, these covariances can be controlled by means of the usual Davydov covariance inequality for mixing processes (see Rio [23, formula 1.12a, page 10] or Bosq [2, formula 1.11, page 22]). We have

$$\forall i \ne j, \quad |Cov(U_i, U_j)| \le C\alpha(|i-j|). \tag{20}$$

To evaluate $s_n^{cov}$, we use the technique developed by Masry [18]. Taking $\varphi_n = \lceil (n^{-1}h_n)^{-1/\nu} \rceil$ (where $\lceil . \rceil$ denotes the smallest integer greater than the argument), we can write

$$s_n^{cov} = \sum_{0 < |i-j| \le \varphi_n} |Cov(U_i, U_j)| + \sum_{|i-j| > \varphi_n} |Cov(U_i, U_j)|. \tag{21}$$

First, applying the upper bound (19) to the first covariance term in (21), we get

$$\sum_{0 < |i-j| \le \varphi_n} |Cov(U_i, U_j)| \le Cnh^2\varphi_n. \tag{22}$$

For the second term, thanks to (20) we get

$$\sum_{|i-j|>\varphi_n} |Cov(U_i, U_j)| \leq C \sum_{|i-j|>\varphi_n} \alpha(|i-j|)$$
$$\leq Cn^2 \alpha(\varphi_n). \tag{23}$$

According to the right-hand side of $(H'2)$, using $(M3)$, (22) and (23), we get

$$s_n^{cov} = O(nh). \tag{24}$$

Finally, (18) and (24) lead directly to $s_n^2 = O(nh)$.

This is enough to study the quantity $Q_{2n}$, since for $\varepsilon$ and $r$ as in (16) and Taylor expansion of $\log(1+x)$ allows us to write that

$$Q_{2n} = C\omega_n^{-1} \exp\left[-\frac{r}{2}\log\left(1 + \frac{\varepsilon_0^2 nh \log n}{rs_n}\right)\right]$$
$$\leq Cn^{\frac{1}{2\beta} - C'\varepsilon_0^2} h^{-(1+\frac{1}{2\beta})}$$
$$= Cn^{\frac{1}{2\beta} - C'\varepsilon_0^2} h^{-\frac{1}{2\beta}(2\beta+1+(\nu+1)\beta)} h^{\frac{\nu+1}{2}}.$$

By using $(H'2)$ and $(M3)$, the later can be made as a general term of a convergent series. Hence $\sum_{n\geq 1}(Q_{1n} + Q_{2n}) < \infty$, and therefore by Borel-Cantelli's Lemma, we have

$$\mathcal{I}_{1n} = O\left(\sqrt{\frac{\log n}{nh}}\right), \quad \mathbf{P}-a.s. \quad \text{as} \quad n \to \infty$$

On the other hand, the bias term $\mathcal{I}_{2n}$ does not depend on the mixing structure. We prove its convergence by using a change of variable and a Taylor expansion (see Lemma 6.1 in Lemdani *et al.* [16]). We get, under $(K3)$ and $(D1)$

$$\mathcal{I}_{2n} = O\left(h^2\right), \quad \mathbf{P}-a.s. \quad \text{as} \quad n \to \infty$$

Hence, replacing $\mathcal{I}_{1n}$ and $\mathcal{I}_{2n}$ in (12), we get the result. ∎

The following Lemma is Lemma 4.2 in Ould-Saïd and Tatachak [22], in which they state a rate of convergence for $\mu_n$ under $\alpha$-mixing hypothesis, which is interesting in itself, similar to that established in the iid case by He and Yang [13].

**Lemma 5.2** *Under assumptions* $(M)$, *we have*

$$|\mu_n - \mu| = O\left(\sqrt{\frac{\log \log n}{n}}\right), \quad \mathbf{P} - a.s. \quad as \quad n \to \infty.$$

**Proof.** See Lemma 4.2 in Ould-Saïd and Tatachak [22]. ∎

Adapting (9), define

$$\tilde{F}_{1,n}(x,y) := \frac{\mu}{nh} \sum_{i=1}^{n} \frac{1}{G(Y_i)} K\left(\frac{x - X_i}{h}\right) H\left(\frac{y - Y_i}{h}\right). \quad (25)$$

**Lemma 5.3** *Under the assumptions of Lemma 5.1 and* $(K2)$, *we have,*

$$\sup_{x \in \Omega} \sup_{a \leq y \leq b} \left|F_{1,n}(x,y) - \tilde{F}_{1,n}(x,y)\right| = O\left(\sqrt{\frac{\log \log n}{n}}\right), \quad \mathbf{P} - a.s. \quad as \quad n \to \infty$$

**Proof.** Under $(K2)$, the df $H$ is bounded by 1. Hence

$$\left|F_{1,n}(x,y) - \tilde{F}_{1,n}(x,y)\right| \leq \left\{\frac{|\mu_n - \mu|}{G_n(a_F)} + \frac{\mu \sup_{y \geq a_F} |G_n(y) - G(y)|}{G_n(a_F)G(a_F)}\right\} |v_n^*(x)|.$$

From Lemma 5.2, we have

$$|\mu_n - \mu| = O\left(\sqrt{\frac{\log \log n}{n}}\right) \quad \mathbf{P} - a.s. \quad as \quad n \to \infty$$

Moreover, $G_n(a_F) \overset{\mathbf{P}-a.s.}{\underset{as}{\to}} \overset{n \to \infty}{} G(a_F) > 0$. In the same way and using Remark 6 in Woodroofe [31] we get

$$\sup_{y \geq a_F} |G_n(y) - G(y)| = O\left(n^{-1/2}\right) \quad \mathbf{P} - a.s. \quad as \quad n \to \infty$$

Combining these last results with Lemma 5.1, we achieve the proof. ∎

**Lemma 5.4** *Under assumptions* $(K)$, $(M)$, $(D3)$, $(D4)$, *and* $(H)$, *we have,*

$$\sup_{x \in \Omega} \sup_{a \leq y \leq b} \left|\tilde{F}_{1,n}(x,y) - \mathbf{E}\left[\tilde{F}_{1,n}(x,y)\right]\right| = O\left(\sqrt{\frac{\log n}{nh}}\right), \quad \mathbf{P}-a.s. \quad as \quad n \to \infty$$

**Proof.** The proof is analogous to that in Lemma 5.1. We give only the leading lines. As $\Omega$ and $[a, b]$ are compact sets, then they can be covered by a finite number $l_n$ and $d_n$ of intervals $I_1, ..., I_{l_n}$ and $J_1, ..., J_{d_n}$ of length $\omega_n$ as in Lemma 6.1 and $\lambda_n = \left(n^{-1}h^{2\beta}\right)^{\frac{1}{2\beta}}$ and centers $x_1, ..., x_{l_n}$ and $y_1, ..., y_{d_n}$ respectively. Since $\Omega$ and $[a, b]$ are bounded, there exist two constant $C_1$ and $C_2$ such that $l_n\omega_n \leq C_1$ and $d_n\lambda_n \leq C_2$. Hence for any $(x, y) \in \Omega \times [a, b]$, there exist $x_k$ and $y_j$ such that $|x - x_k| \leq \omega_n$ and $|y - y_j| \leq \lambda_n$. Thus we have the following decomposition

$$\sup_{x \in \Omega} \sup_{y \in [a,\, b]} \left|\tilde{F}_{1,n}(x, y) - \mathbf{E}\left[\tilde{F}_{1,n}(x, y)\right]\right|$$

$$\leq \max_{1 \leq k \leq l_n} \sup_{x \in I_k} \sup_{y} \left|\tilde{F}_{1,n}(x, y) - \tilde{F}_{1,n}(x_k, y)\right|$$

$$+ \max_{1 \leq k \leq l_n} \max_{1 \leq j \leq d_n} \sup_{y \in J_j} \left|\tilde{F}_{1,n}(x_k, y) - \tilde{F}_{1,n}(x_k, y_j)\right|$$

$$+ \max_{1 \leq k \leq l_n} \max_{1 \leq j \leq d_n} \left|\tilde{F}_{1,n}(x_k, y_j) - \mathbf{E}\left[\tilde{F}_{1,n}(x_k, y_j)\right]\right|$$

$$+ \max_{1 \leq k \leq l_n} \max_{1 \leq j \leq d_n} \sup_{y \in J_j} \left|\mathbf{E}\left[\tilde{F}_{1,n}(x_k, y_j)\right] - \mathbf{E}\left[\tilde{F}_{1,n}(x_k, y)\right]\right|$$

$$+ \max_{1 \leq k \leq l_n} \sup_{x \in I_k} \sup_{y} \left|\mathbf{E}\left[\tilde{F}_{1,n}(x_k, y)\right] - \mathbf{E}\left[\tilde{F}_{1,n}(x, y)\right]\right|$$

$$=: \mathcal{J}_{1n} + \mathcal{J}_{2n} + \mathcal{J}_{3n} + \mathcal{J}_{4n} + \mathcal{J}_{5n}.$$

Firstly, concerning $\mathcal{J}_{1n}$ and $\mathcal{J}_{5n}$, assumptions $(K1)$ and $(K2)$ yield

$$\sup_{x \in I_k} \sup_{y} \left|\tilde{F}_{1,n}(x, y) - \tilde{F}_{1,n}(x_k, y)\right| \leq \frac{C\mu\omega_n^\beta}{G(a_F)h^{1+\beta}} \sup_{y} \left|H\left(\frac{y - Y_i}{h}\right)\right| = O\left((nh)^{-1/2}\right).$$

Hence, by $(H1)$ we get

$$\sqrt{\frac{nh_n}{\log n}} \sup_{x \in \Omega} \sup_{y \in [a,\, b]} \left|\tilde{F}_{1,n}(x, y) - \tilde{F}_{1,n}(x_k, y)\right| = o(1). \tag{26}$$

Similarly, we obtain for $\mathcal{J}_{2n}$ and $\mathcal{J}_{4n}$

$$\sup_{y \in J_j} \left|\tilde{F}_{1,n}(x_k, y) - \tilde{F}_{1,n}(x_k, y_j)\right| \leq \frac{C\mu\lambda_n^\beta}{G(a_F)h^{1+\beta}} \left|K\left(\frac{x_k - X_i}{h}\right)\right| = O\left(\left(nh_n^2\right)^{-1/2}\right).$$

Again, by $(H1)$ we get

$$\sqrt{\frac{nh}{\log n}} \sup_{x \in \Omega} \sup_{y \in [a,\, b]} \left|\tilde{F}_{1,n}(x_k, y) - \tilde{F}_{1,n}(x_k, y_j)\right| = o(1). \tag{27}$$

As to $\mathcal{J}_{3n}$, for all $\varepsilon > 0$ we have

$$\mathbf{P}\left\{\max_{1\leq k\leq l_n}\max_{1\leq j\leq d_n}\left|\tilde{F}_{1,n}(x_k,y_j) - \mathbf{E}\left[\tilde{F}_{1,n}(x_k,y_j)\right]\right| > \varepsilon\right\}$$
$$\leq l_n d_n \mathbf{P}\left\{\left|\tilde{F}_{1,n}(x_k,y_j) - \mathbf{E}\left[\tilde{F}_{1,n}(x_k,y_j)\right]\right| > \varepsilon\right\}. \tag{28}$$

Set, for any $i \geq 1$,

$$\Psi_i(x_k,y_j) := \frac{\mu}{nh}\left\{\frac{1}{G(Y_i)}K\left(\frac{x_k-X_i}{h}\right)H\left(\frac{y_j-Y_i}{h}\right)\right.$$
$$\left. - \mathbf{E}\left[\frac{1}{G(Y_i)}K\left(\frac{x_k-X_1}{h}\right)H\left(\frac{y_j-Y_i}{h}\right)\right]\right\}.$$

Under $(K1)$ and $(K2)$, the rv's $V_i := nh\Psi_i(x_k, y_j)$ are centered and bounded by $\frac{2\mu M_0 M_1}{G(a_F)} =: C < \infty$. Then, applying again Fuk-Nagaev inequality, we obtain that, for all $\varepsilon > 0$ and $r > 1$,

$$\mathbf{P}\left\{\max_{1\leq k\leq l_n}\max_{1\leq j\leq d_n}\left|\sum_{i=1}^n \Psi_i(x_k,y_j)\right| > \varepsilon\right\}$$
$$= \mathbf{P}\left\{\max_{1\leq k\leq l_n}\max_{1\leq j\leq d_n}\left|\sum_{i=1}^n V_i\right| > nh\varepsilon\right\}$$
$$\leq C_1 C_2 (\omega_n \lambda_n)^{-1}\left\{\frac{n}{r}\left(\frac{2r}{\varepsilon nh}\right)^{\nu+1} + \left(1 + \frac{\varepsilon^2 n^2 h^2}{rs_n^2}\right)^{-\frac{r}{2}}\right\}$$
$$\leq Cn^{\frac{1}{\beta}}h_n^{-\left(\frac{1}{2\beta}+2\right)}\frac{n}{r}\left(\frac{r}{\varepsilon nh}\right)^{\nu+1} + Cn^{\frac{1}{\beta}}h^{-\left(\frac{1}{2\beta}+2\right)}\left(1 + \frac{\varepsilon^2 n^2 h^2}{rs_n^2}\right)^{-\frac{r}{2}}$$
$$=: \mathcal{J}_{31n} + \mathcal{J}_{32n}, \tag{29}$$

where $s_n^2 = \sum_{1\leq i\leq n}\sum_{1\leq j\leq n}|Cov(V_i,V_j)|$.

By taking $\varepsilon$ and $r$ as in (16), we get

$$\mathcal{J}_{31n} = C\varepsilon_0^{-(\nu+1)}n^{1+\frac{1}{\beta}-\frac{\nu+1}{2}}(\log n)^{\nu(1+\delta)-\frac{\nu+1}{2}}h^{\frac{-1}{2\beta}(1+4\beta+(\nu+1)\beta)}.$$

Then, using $(H1)$ and $(H2)$ we get

$$\mathcal{J}_{31n} \leq C\varepsilon_0^{-(\nu+1)}n^{-1-\frac{\eta}{2\beta}\left(1+4\beta+(\nu+1)\beta-\frac{2}{\eta}\right)}(\log n)^{\nu(1+\delta)-\frac{\nu-1}{2}}.$$

Hence, the condition upon $\beta$ and for any $\eta$ as in $(H2)$, $\mathcal{J}_{31n}$ is the general term of a finite-sum series.

Let us now examine the term $\mathcal{J}_{32n}$. First, we have to calculate

$$s_n^2 = nVar(V_1) + \sum_{i \neq j} |Cov(V_i, V_j)|.$$

We have

$$Var(V_1) = \mathbf{E}\left[\frac{\mu^2}{G^2(Y_1)} K^2\left(\frac{x_k - X_i}{h}\right) H^2\left(\frac{y_j - Y_1}{h}\right)\right]$$
$$- \mathbf{E}^2\left[\frac{\mu}{G(Y_i)} K\left(\frac{x_k - X_i}{h}\right) H\left(\frac{y_j - Y_1}{h}\right)\right]$$
$$=: \mathcal{V}_1 + \mathcal{V}_2.$$

Remark that

$$\mathbf{E}\left[\frac{\mu^2}{G^2(Y_1)} H^2\left(\frac{y_j - Y_1}{h}\right)\bigg| X_1\right] = \int \frac{\mu^2}{G^2(y_1)} H^2\left(\frac{y_j - y_1}{h}\right) f^*(y_1|X_1) dy_1$$
$$= \int \frac{\mu}{G(y_1)} H^2\left(\frac{y_j - y_1}{h}\right) f(y_1|X_1) dy_1$$
$$= \mathbf{E}\left[\frac{\mu}{G(Y_1)} H^2\left(\frac{y_j - Y_1}{h}\right)\right].$$

Then

$$\mathcal{V}_1 = \mathbf{E}\left[\frac{\mu}{G(Y_1)} K^2\left(\frac{x_k - X_i}{h}\right) H^2\left(\frac{y_j - Y_1}{h}\right)\right]$$
$$\leq \frac{\mu}{G(a_F)} \mathbf{E}\left[K^2\left(\frac{x_k - X_i}{h}\right)\right]$$
$$\leq \frac{\mu h_n}{G(a_F)} \int K^2(r) v^*(x_k - rh) dr.$$

Under $(K3:ii)$ and $(D1)$, we have $\mathcal{V}_1 = O(h)$. An analogous developpement gives that $\mathcal{V}_2 = O(h^2)$, which implies $Var(V_1) = O(h)$.

On the one hand, $(M1)$, $(K1)$ and $(K2)$ lead to

$$\begin{aligned}|Cov(V_i,V_j)| = &\left|\int\int\int\int \frac{\mu}{G(r)}K\left(\frac{x_k-u}{h}\right)H\left(\frac{y_j-r}{h}\right)\frac{\mu}{G(t)}K\left(\frac{x_k-s}{h}\right)\right.\\&\times H\left(\frac{y_j-t}{h}\right)f^*_{1,1,j+1,j+1}(u,r,s,t)dudrdsdt\\&-\int\int \frac{\mu}{G(r)}K\left(\frac{x_k-u}{h}\right)H\left(\frac{y_j-r}{h}\right)f^*(u,r)dudr\\&\times\left.\int\int \frac{\mu}{G(t)}K\left(\frac{x_k-s}{h}\right)H\left(\frac{y_j-t}{h}\right)f^*(s,t)dsdt\right|\\\leq &\frac{\mu^2}{G^2(a_F)}\int\int\int\int \left|K\left(\frac{x_k-u}{h}\right)H\left(\frac{y_j-r}{h}\right)K\left(\frac{x_k-s}{h}\right)H\left(\frac{y_j-t}{h}\right)\right.\\&\times \left.\left(f^*_{1,1,j+1,j+1}(u,r,s,t)-f^*(u,r)f^*(s,t)\right)\right|dudrdsdt.\end{aligned}$$

Using Assumption $(D3)$ and by a change of variable, it follows that

$$|Cov(V_i,V_j)| = O\left(h^4\right). \tag{30}$$

On the other hand, from a result in Bosq [2, p.22], we have

$$|Cov(V_i,V_j)| = O\left(\alpha\left(|i-j|\right)\right). \tag{31}$$

Then to evaluate $\sum_{i\neq j}|Cov(V_i,V_j)|$, the idea is to introduce a sequence of integers $\varphi_n$ the same as in Lemma 5.1, and using (30) for the nearest and (31) for the farest integer $i$ and $j$. Then we get

$$\begin{aligned}\sum_{i\neq j}|Cov(V_i,V_j)| &= \sum\sum_{0<|i-j|\leq\varphi_n}|Cov(V_i,V_j)| + \sum\sum_{|i-j|>\varphi_n}|Cov(V_i,V_j)|\\&\leq \sum\sum_{0<|i-j|\leq\varphi_n}h^4 + \sum\sum_{|i-j|>\varphi_n}\alpha\left(|i-j|\right)\\&\leq Cn\varphi_n h^4 + Cn^2\alpha\left(\varphi_n\right).\end{aligned}$$

The right-hand side of $(H2)$ and $(M3)$, one has $\sum_{i\neq j}|Cov(V_i,V_j)| = O(nh)$.
So $s_n = O(nh)$.

Consequently, by taking $r$ and $\varepsilon$ as in (16) and using Taylor expansion of $\log(1+x)$, the term $\mathcal{J}_{32n}$ becomes

$$\begin{aligned}\mathcal{J}_{32n} &\leq Cn^{\frac{1}{\beta}}h^{-\left(\frac{1}{2\beta}+2\right)}\exp\left\{-\frac{1}{2}\varepsilon_0^2\log n\right\}\\&= Cn^{\frac{1}{\beta}-C\varepsilon_0^2}h^{-\left(\frac{1}{2\beta}+2\right)}.\end{aligned}$$

By using $(H2)$ and $(M3)$, the later can be made as a general term of summable series. Thus $\sum_{n\geq 1}(\mathcal{J}_{31n} + \mathcal{J}_{32n}) < \infty$. Then by Borel-Cantelli's Lemma, the first term of $(29)$ goes to zero $a.s.$ and for $n$ large enough, we have $\mathcal{J}_{3n} = O\left(\sqrt{\frac{\log n}{nh}}\right)$, this complet the proof of the Lemma. ∎

**Lemma 5.5** *Under assumptions $(K3)$ and $(D4)$ we have,*

$$\sup_{x\in\Omega}\sup_{a\leq y\leq b}\left|\mathbf{E}\left[\tilde{F}_{1,n}(x,y)\right] - F_1(x,y)\right| = O\left(h^2\right), \quad \mathbf{P}-a.s. \quad as \quad n\to\infty$$

**Proof.** The bias terms do not depend on the mixing structure. The proof of Lemma 5.5 is similar to that of Lemma 6.2 in Lemdani *et al.* [16], hence its proof is omitted. ∎

The next Lemma gives the uniform convergence with rate of the estimator $v_n(x)$ defined in (7).

**Lemma 5.6** *Under the assumptions of Lemma 5.1 and condition $(D5)$, we have*

$$\sup_{x\in\Omega}|v_n(x) - v(x)| = O\left(\max\left\{\sqrt{\frac{\log n}{nh}}, h^2\right\}\right), \quad \mathbf{P}-a.s. \quad as \quad \to\infty.$$

**Proof.** Adapting (7), define

$$\tilde{v}_n(x) = \frac{\mu}{nh}\sum_{i=1}^{n}\frac{1}{G(Y_i)}K\left(\frac{x-X_i}{h}\right). \qquad (32)$$

We have

$$\sup_{x\in\Omega}|v_n(x) - v(x)| \leq \sup_{x\in\Omega}|v_n(x) - \tilde{v}_n(x)|$$
$$+ \sup_{x\in\Omega}|\tilde{v}_n(x) - \mathbf{E}\left[\tilde{v}_n(x)\right]|$$
$$+ \sup_{x\in\Omega}|\mathbf{E}\left[\tilde{v}_n(x)\right] - v(x)|$$
$$=: \mathcal{L}_{1n} + \mathcal{L}_{2n} + \mathcal{L}_{3n}.$$

For the first term, using analogous framework as in Lemma 5.3, we get

$$\mathcal{L}_{1n} = O\left(\sqrt{\frac{\log \log n}{n}}\right), \quad \mathbf{P}-a.s. \quad \text{as} \quad \to \infty. \tag{33}$$

In addition, by using the same approach as for $\mathcal{I}_{1n}$ in the proof of Lemma 5.1, we can show that, for $n$ large enough

$$\mathcal{L}_{2n} = O\left(\sqrt{\frac{\log n}{nh}}\right), \quad \mathbf{P}-a.s. \quad \text{as} \quad \to \infty. \tag{34}$$

Finally, a change of variable and a Taylor expansion we get, under $(K3)$ and $(D5)$

$$\mathbf{E}\left[\tilde{v}_n(x)\right] - v(x) = \mathbf{E}\left[\frac{\mu}{nh}\sum_{i=1}^{n}\frac{1}{G(Y_1)}K\left(\frac{x-X_1}{h}\right)\right] - v(x)$$

$$= \frac{1}{h}\int K\left(\frac{x-u}{h}\right)v(u)du - v(x)$$

$$= \frac{h^2}{2}\int r^2 K(r) v''(\tilde{x})dr$$

with $\tilde{x} \in [x - rh, x]$, which yields that

$$\mathcal{L}_{3n} = O\left(h^2\right), \quad \mathbf{P}-a.s. \quad \text{as} \quad n \to \infty. \tag{35}$$

Combining (33), (34) and (35) permit to conclude the proof. ∎

**Proof of Proposition 3.1** In view of (8), we have the following classical decomposition

$$\sup_{x \in \Omega}\sup_{a \leq y \leq b}|F_n(y|x) - F(y|x)| \leq \frac{1}{\gamma - \sup_{x \in \Omega}|v_n(x) - v(x)|}\left\{\sup_{x \in \Omega}\sup_{a \leq y \leq b}|F_{1,n}(x,y) - F_1(x,y)|\right.$$

$$\left. + \gamma^{-1}\sup_{x \in \Omega}\sup_{a \leq y \leq b}|F(y|x)|\sup_{x \in \Omega}|v_n(x) - v(x)|\right\}.$$

Furthermore, we have

$$\sup_{x \in \Omega}\sup_{a \leq y \leq b}|F_{1,n}(x,y) - F_1(x,y)| \leq \sup_{x \in \Omega}\sup_{a \leq y \leq b}\left|F_{1,n}(x,y) - \tilde{F}_{1,n}(x,y)\right|$$

$$+ \sup_{x \in \Omega}\sup_{a \leq y \leq b}\left|\tilde{F}_{1,n}(x,y) - \mathbf{E}\left[\tilde{F}_{1,n}(x,y)\right]\right|$$

$$+ \sup_{x \in \Omega}\sup_{a \leq y \leq b}\left|\mathbf{E}\left[\tilde{F}_{1,n}(x,y)\right] - F_1(x,y)\right|.$$

In conjunction with Lemmas 5.3–5.6, we conclude the proof. ∎

We now embark on the proof of Theorem 3.1.

**Proof of Theorem 3.1** Let $x \in \Omega$. As $F_n(\cdot|x)$ and $F(\cdot|x)$ are continuous, we have $F(q_p(x)|x) = F_n(q_{p,n}(x)|x) = p$. Then

$$\begin{aligned}
|F(q_{p,n}(x)|x) - F(q_p(x)|x)| &\leq |F(q_{p,n}(x)|x) - F_n(q_{p,n}(x)|x)| \\
&\quad + |F_n(q_{p,n}(x)|x) - F(q_p(x)|x)| \\
&\leq |F(q_{p,n}(x)|x) - F_n(q_{p,n}(x)|x)| \\
&\leq \sup_{a \leq y \leq b} |F_n(y|x) - F(y|x)|. \quad (36)
\end{aligned}$$

The consistency of $q_{p,n}(x)$ follows then immediately from Proposition 3.1 in conjunction with the inequality

$$\sum_n \left\{ \sup_{x \in \Omega} |q_{p,n}(x) - q_p(x)| \geq \varepsilon \right\} \leq \sum_n \left\{ \sup_{x \in \Omega} \sup_{a \leq y \leq b} |F_n(y|x) - F(y|x)| \geq \beta \right\}.$$

For the second part, a Taylor expansion of $F(.|.)$ in neighborhood of $q_p$, implies that

$$F(q_{p,n}(x)|x) - F(q_p(x)|x) = (q_{p,n}(x) - q_p(x)) f(\tilde{q}_p(x)|x) \quad (37)$$

where $\tilde{q}_p$ is between $q_p$ and $q_{p,n}$ and $f(.|x)$ is the conditional density of $Y$ given $X = x$. Then, from the behavior of $F(q_{p,n}(x)|x) - F(q_p(x)|x)$ as $n$ goes to infinity, it is easy to obtain asymptotic results for the sequence $(q_{p,n}(x) - q_p(x))$. By (37) we have

$$\sup_{x \in \Omega} |q_{p,n}(x) - q_p(x)| |f(\tilde{q}_p(x)|x)| \leq \sup_{x \in \Omega} \sup_{a \leq y \leq b} |F_n(y|x) - F(y|x)|.$$

The result follows from $(D4)$ and the Proposition 3.1. Here we point out that, if $f(\tilde{q}_p(x)|x) = 0$, for some $x \in \Omega$, we should increase the order of Taylor expansion to obtain the consistency of $q_{p,n}(x)$ (with an adapted rate). ∎

**References**

[1] Boente G, Fraiman R (1995) Asymptotic distribution of smoothers based on local means and local medians under dependence. J of Multivariate Anal 54:77–90.

[2] Bosq D (1998) Nonparametric Statistics for Stochastic Processes: Estimation and Prediction. (Second Edition). Lecture Notes in Statistics, 110, Springer Verlag, New York.

[3] Cai, Z. (1998) Kernel density and hazard rate estimation for censored dependent data. J of Multivariate Anal 67:23–34.

[4] Cai, Z (2001) Estimating a distribution function for censored time series data. J Multivariate Anal 78:299-318.

[5] Chaudhuri, P., Doksum, K., Samarov, A (1997) On average derivative quantile regression. Ann Statist, 25:715-744.

[6] Doukhan, P (1994) Mixing: Properties and examples. Lecture Notes in Statistics 85, Springer-Verlag, New York.

[7] Feigelson, E.D., Babu, G.J (1992) Statistical Challenges in Modern Astronomy. Berlin Heidelberg, Springer Verlag, New York.

[8] Ferraty, F., Vieu, P (2006) Nonparametric Functional Data Analysis, Theory and Practice. Springer-Verlag, New York.

[9] Gannoun A (1989) Estimation de la médiane conditionnelle. Thèse de doctorat de l'Université de Paris VI.

[10] Gannoun A (1990) Estimation non paramétrique de la médiane conditionnelle médianogramme et méthode du noyau. Publication de l'Institut de Statistique de l'université de Paris 11-22.

[11] Gürler U, Stute W, Wang, JL (1993) Weak and strong quantile representations for randomly truncated data with applications. Statist Probab Lett, 17:139–148.

[12] He S, Yang, G (1994) Estimating a lifetime distribution under different sampling plan. In S.S. Gupta J.O. Berger (Eds.) Statistical decision theory and related topics 5:73-85 Berlin Heidelberg, Springer Verlag, New York.

[13] He S, Yang, G (1998) Estimation of the truncation probability in the random truncation model. Ann Statist 26:1011–1027.

[14] Lemdani M, Ould-Saïd E (2007) Asymptotic behaviorof the hazard rate kernel estimator under truncated and censored data. Comm. in Statist. Theory & Methods 37-155-174.

[15] Lemdani M, Ould-Saïd E, Poulin, N (2005) Strong representation of the

quantile function for left-truncated and dependent data. Math Meth Statist 14:332–345

[16] Lemdani M, Ould-Saïd E, Poulin, N (2008) Prediction for a left truncated model via estimation of the conditional quantile. (In press in J. Multivariate Anal.).

[17] Lynden-Bell D (1971) A method of allowing for known observational selection in small samples applied to 3CR quasars. *Monthly Notices Roy Astron Soc* 155:95–118.

[18] Masry E (1986) Recursive probability density estimation for weakly dependent processes. IEEE Transactions on Information Theory 32:254-267.

[19] Masry E, Tj$\phi$stheim, D (1995) Nonparametric estimation and identification of nonlinear time series. Econometric Theory 11:258-289.

[20] Mehra KL, Rao MS, Upadrasta, SP (1991) A smooth conditional quantile estimator and related applications of conditional empirical processes. J Multivariate Anal 37:151–179.

[21] Ould-Saïd E, Lemdani, M (2006) Asymptotic properties of a nonparametric regression function estimator with randomly truncated data. Ann Instit Statist Math 58:357–378.

[22] Ould-Saïd E, Tatachak, A (2007) Strong uniform consistency rate for the kernel mode under strong mixing hypothesis and left-truncation. (In press in Comm. Statist. Theory & Method).

[23] Rio E (2000) Théorie Asymptotique des Processus Aléatoires Faiblement Dépendants (In French). Mathématiques et Applications, 31, Springer Verlag, New York.

[24] Roussas GG (1969). Nonparametric estimation of the transition distribution function of a Markov process. Ann Math Statist 40:1386–140.

[25] Samanta M (1989) Nonparametric estimation of conditional quantiles. *Statist Probab Lett* **7**, 407–412.

[26] Schlee W (1982) Estimation non paramétrique du $\alpha-$quantile conditionnel. Statistique et Analyse des données 1:32–47.

[27] Stone C (1977) Consistent nonparametric regression. Ann Statist 5:595–645.

[28] Stute W (1993) Almost sure representation of the product-limit estimator for truncated data. Ann. Statits. 21:146-156.

[29] Tsai WY, Jewell NP, Wang, MC (1987). A note on the product-limit estimator under right censoring and left truncation. Biometrika 74:883-886.

[30] Wang MC, Jewell NP, Tsai, WY (1986) Asymptotic properties of the product-limit estimate under random truncation. Ann Statist 14:1597-1605.
[31] Woodroofe M (1985) Estimating a distribution function with truncated data. Ann Statist 13:163-177.
[32] Xiang X (1996) A kernel estimator of a conditional quantile. J Multivariate Anal 59:206–216.